\documentclass[letterpaper, 10pt, conference]{ieeeconf}      % Use this line for a4
\IEEEoverridecommandlockouts                              % This command is only
\usepackage{amsmath,graphicx,amsfonts,amssymb,epsfig,mathrsfs,enumerate}
\usepackage{subcaption,color}
\usepackage{algorithm,algpseudocode,algorithmicx}
\usepackage{cite,url,framed,bm,dsfont}

\newtheorem{theorem}{Theorem}

\newtheorem{lemma}{Lemma}
\newtheorem{proposition}{Proposition}
\newtheorem{remark}{Remark}

\usepackage{tikz}
\usetikzlibrary{calc,trees,positioning,arrows,chains,shapes.geometric,decorations.pathreplacing,decorations.pathmorphing,shapes, matrix,shapes.symbols}

\newcommand{\bw}{{\bm{w}}}
\newcommand{\bx}{{\bm{x}}}

\newcommand{\bbf}{{\bm{f}}}

\newcommand{\bu}{{\bm{u}}}

\newcommand{\bn}{{\bm{n}}}

\renewcommand{\dagger}{{T}}
\newcommand{\differential}{{\rm{d}}}
\newcommand{\uopt}{\bm{u}^{\rm{opt}}}
\newcommand{\ropt}{\rho^{\rm{opt}}}
\newcommand{\tx}{\Tilde{x}}
\newcommand{\ty}{\Tilde{y}}

\newcommand*{\qed}{\hfill\ensuremath{\blacksquare}}

\allowdisplaybreaks

\title{\LARGE\textbf{Reflected Schr\"{o}dinger Bridge: Density Control with Path Constraints
}
}
\author{Kenneth F. Caluya, and Abhishek Halder% <-this % stops a space
\thanks{Kenneth F. Caluya, and Abhishek Halder are with the Department of Applied Mathematics, University of California, Santa Cruz, CA 95064, USA,
        {\tt\small{\{kcaluya,ahalder\}@ucsc.edu}}. This research was partially supported by NSF award 1923278.%
}}

\begin{document}

\maketitle
\thispagestyle{empty}
\pagestyle{empty}

\def\spacingset#1{\def\baselinestretch{#1}\small\normalsize}
\spacingset{1}

\title{\huge{Reflecting Schrodinger Bridge}}
\begin{abstract}
How to steer a given joint state probability density function to another over finite horizon subject to a controlled stochastic dynamics with hard state (sample path) constraints? In applications, state constraints may encode safety requirements such as obstacle avoidance. In this paper, we perform the feedback synthesis for minimum control effort density steering (a.k.a. Schr\"{o}dinger bridge) problem subject to state constraints. We extend the theory of Schr\"{o}dinger bridges to account the reflecting boundary conditions for the sample paths, and provide a computational framework building on our previous work on proximal recursions, to solve the same.
\end{abstract}

\section{Introduction}
We consider finite horizon feedback steering of an ensemble of trajectories subject to a controlled stochastic differential equation (SDE) with endpoint joint state probability density function (PDF) constraints -- a topic of growing interest in the systems-control literature. Motivating applications include belief space motion planning for vehicular autonomy, and the steering of robotic or biological swarms via decentralized feedback. While early contributions focused on the covariance control \cite{hotz1987covariance,skelton1989covariance,grigoriadis1997minimum}, more recent papers \cite{chen2015optimal1,chen2015optimal2,chen2016optimal} addressed the optimal feedback synthesis for steering an arbitrary prescribed initial joint state PDF to another prescribed terminal joint state PDF subject to controlled linear dynamics, and revealed the connections between the associated stochastic optimal control problem, the theory of optimal mass transport \cite{villani2003topics}, and the Schr\"{o}dinger bridge \cite{schrodinger1931umkehrung,schrodinger1932theorie}. Follow up works have accounted terminal cost \cite{halder2016finite}, input constraints \cite{bakolas2018finite,okamoto2019input}, output feedback \cite{bakolas2017covariance}, and some nonlinear dynamics \cite{caluya2019finite,caluya2019finiteMulti,caluya2019wasserstein}. As for the state or path constraints, prior work \cite{okamoto2018optimal} incorporated the same in soft probabilistic sense. The contribution of the present paper is to account \emph{hard deterministic path constraints} in the problem of minimum effort finite horizon PDF steering via feedback synthesis. This can be intuitively phrased as the ``hard safety with soft endpoint" problem.

The main idea underlying the ensuing development is to modify the unconstrained It\^{o} SDEs to the ``reflected It\^{o} SDEs" \cite{ikeda1961construction,watanabe1971stochastic,lions1984stochastic,harrison1987multidimensional}, i.e., the controlled sample paths in the state space (in addition to the control-affine deterministic drift) are driven by two stochastic processes: a Wiener process, and a local time stochastic process. The latter enforces the sample paths in the state space to satisfy the deterministic non-strict\footnote{There is no loss of generality in allowing the sample paths to satisfy \emph{non-strict} path containment in given $\mathcal{X} \subset \mathbb{R}^n$ since \emph{strict} containment can be enforced by reflecting them from $\epsilon$-inner boundary layer of $\partial\mathcal{X}$ for $\epsilon$ small enough.} path containment constraints at all times. These considerations engender a Schr\"{o}dinger bridge-like formulation--referred hereafter as the \emph{Reflected Schr\"{o}dinger Bridge Problem (RSBP)}--which unlike its classical counterpart, has extra boundary conditions involving the gradients of the so-called Schr\"{o}dinger factors. We show how recent developments in contraction mapping w.r.t. the Hilbert metric, and the proximal recursion over the Schr\"{o}dinger factors can be harnessed to solve the RSBP.

\section{Reflected Schr\"{o}dinger Bridge Problem}\label{SecRSBP}

\subsection{Formulation}\label{SubsecFormulation}
Consider a connected, smooth\footnote{More precisely, there exists $\xi\in C_{b}^{2}\left(\mathbb{R}^{n}\right)$ such that $\mathcal{X} \equiv \{\bx\in\mathbb{R}^{n}\mid \xi(\bx) > 0\}$ with boundary $\partial\mathcal{X}\equiv\{\bx\in\mathbb{R}^{n}\mid \xi(\bx) = 0\}$.} and bounded domain $\mathcal{X} \subset \mathbb{R}^n$. Let $\overline{\mathcal{X}}:=\mathcal{X}\cup\partial\mathcal{X}$ denote the closure of $\mathcal{X}$. For time $t\in[0,1]$, consider the stochastic control problem
 \begin{subequations} \label{ReflectedOCPM}  
\begin{align} 
& \underset{\bu \in \mathcal{U}} {\text{inf}}
& &  \mathbb{E} \left\{\int_{0}^{1} \frac{1}{2} \| \bu(t,\bx_{t}^{\bu})  \|_{2}^{2}  \: \differential t \right\}\label{ReflectedOCPobjM}  \\
& \text{subject to}
&  & \differential \bx_t^{\bu} = \bbf(t,\bx_t^{\bu}) \:  \differential t + \bu(t,\bx_t^{\bu}) \: \differential t \nonumber \\
& &  & \quad \qquad +\sqrt{2\theta} \: \differential \bw_t +\bm{n}(\bx_t^{\bu}) \differential \gamma_t, \label{ReflectedOCPdynM}  \\  
& & & \!\!\!\!\!\!\!\!\!\!\!\!\!\!\!\!\!\!\!\!\!\!\!\!\bx_{0}^{\bu} :=  \bx_t^{u}\left(t=0\right) \sim  \rho_{0}, \quad \bx_1^{\bu}:=  \bx_t^{u}\left(t=1\right) \sim \rho_{1}, \label{ReflectedOCPbcM}  
\end{align}
\end{subequations}
where $\bw_t$ is the standard Wiener process in $\mathbb{R}^{n}$, the controlled state $\bx_{t}^{\bu}\in\overline{\mathcal{X}}$, and the endpoint joint state PDFs $\rho_{0},\rho_{1}$ are prescribed\footnote{The notation $\bx \sim \rho$ means that the random vector $\bx$ has joint PDF $\rho$.} such that their supports are in $\overline{\mathcal{X}}$, both are everywhere nonnegative, have finite second moments, and $\int\rho_{0}=\int\rho_{1}=1$. The parameter $\theta > 0$ is referred to as the thermodynamic temperature, and the expectation operator $\mathbb{E}\{\cdot\}$ in (\ref{ReflectedOCPobjM}) is w.r.t. the law of the controlled state $\bx_{t}^{\bu}$. The set $\mathcal{U}$ consists of all admissible feedback policies $\bm{u}(t,\bx_t^{u})$, given by $\mathcal{U}:=\{\bm{u}: [0,1]\times\overline{\mathcal{X}} \mapsto \mathbb{R}^{n} \mid \|\bm{u}\|_{2}^{2} < \infty, \bu(t,\cdot)\in\text{Lipschitz}\left(\overline{\mathcal{X}}\right)\:\text{for all}\:t\in[0,1]\}$. We assume that the prior drift vector field $\bm{f}$ is bounded Borel measurable in $(t,\bx_{t}^{\bu})\in[0,1]\times\overline{\mathcal{X}}$, and Lipschitz continuous w.r.t. $\bx_{t}^{\bu}\in\overline{\mathcal{X}}$. The vector field $\bm{n}$ is set to be the inward unit normal to the boundary $\partial\mathcal{X}$, and gives the direction of reflection. Furthermore, for $t\in[0,1]$, $\gamma_t$ is \emph{minimal local time}: a continuous, non-negative and non-decreasing stochastic process \cite{glynn2018rate,harrison2013brownian,pilipenko2014introduction} that restricts $\bx_t^{u}$ to the domain $\mathcal{X}$, with $\gamma_{0}\equiv 0$. Specifically, letting $\mathds{1}_{\{\}}$ denote the indicator function of the subscripted set, we have
\begin{equation}
\gamma_t = \int_{0}^{t}\mathds{1}_{\{\bx_s^{\bu} \in \partial\mathcal{X}\}}\: \differential \gamma_s, \quad \int_{0}^1 \mathds{1}_{\{\bx_t^{\bu} \notin \partial \mathcal{X}\}}\: \differential \gamma_t =0,
\label{gammat}
\end{equation}
which is to say that the process $\gamma_t$ only increases at times $t\in[0,1]$ when $\bx_t^{\bu}$ hits the boundary, i.e.,  when $\bx_t^{\bu} \in  \partial \mathcal{X} $. Thus, (\ref{ReflectedOCPdynM}) is a controlled reflected SDE, and the tuple $(\bx_t^{\bu},\gamma_t)$ solves the \emph{Skorokhod problem} \cite{skorokhod1961stochastic,skorokhod1962stochastic,kruk2007explicit}.
We point the readers to \cite{lions1984stochastic} for the proof of existence and uniqueness of solutions to (\ref{ReflectedOCPdynM}) under the stated regularity assumptions.

To formalize the probabilistic setting of the problem at hand, let $\Omega$ be the space of continuous functions $\omega: [0,1] \mapsto \overline{\mathcal{X}}$. We view $\Omega$ as a complete separable metric space endowed with the topology of uniform convergence on compact time intervals. With $\Omega$, we associate the $\sigma$-algebra $\mathscr{F} = \sigma\{\omega(s) \mid 0 \leq s \leq 1\}$. Consider the complete filtered probability space $\left(\Omega,\mathscr{F},\mathbb{P}\right)$ with filtration $\mathscr{F}_{t} = \sigma\{\omega(s)\mid 0\leq s \leq t \leq 1\}$ wherein ``complete" means that $\mathscr{F}_{0}$ contains all $\mathbb{P}$-null sets, and $\mathscr{F}_{t}$ is right continuous. The processes $\bw_{t}$, $\bx_{t}^{\bu}$ (for a given feedback policy $\bu$) and $\gamma_{t}$ are $\mathscr{F}_{t}$-adapted (i.e., non-anticipating) for $t\in[0,1]$. In (\ref{ReflectedOCPbcM}), the random vectors $\bx_{0}^{\bu}$ and $\bx_{1}^{\bu}$ are respectively $\mathscr{F}_{0}$-measurable and $\mathscr{F}_{1}$-measurable. 

Denote the Euclidean gradient operator as $\nabla$, the inner product as $\langle\cdot,\cdot\rangle$, and the Laplacian as $\Delta$. Letting 
\[\mathcal{L} := \theta\Delta + \langle \bbf + \bu,\nabla\rangle,\]
the law of the sample path of (\ref{ReflectedOCPdynM}) can be characterized \cite{stroock1971diffusion} as follows: for each $\bx\in\overline{\mathcal{X}}$, there is a \emph{unique} probability measure $\mathbb{P}_{\bx}^{\mu}$ on $\Omega$ such that (i) $\mathbb{P}_{\bx}^{\mu}\left(\bx_t^{\bu}(t=0) = \bx\right) = 1$, (ii) for any $\phi\in C_{c}^{1,2}\left([0,1];\overline{\mathcal{X}}\right)$ whose inner normal derivative on $\partial\mathcal{X}$ is nonnegative,
\[\phi(t,\bx_{t}^{\bu}) - \displaystyle\int_{0}^{t}\left(\dfrac{\partial\phi}{\partial s} + \mathcal{L}\phi\right)\left(s,\bx_{s}^{\bu}\right)\:\differential s\]
is $\mathbb{P}^{\mu}_{\bx}$-submartingale, and (iii) there is a continuous, nonnegative, nondecreasing stochastic process $\gamma_{t}$ satisfying (\ref{gammat}). As a consequence \cite[p. 196]{stroock1971diffusion} of this characterization it follows that the process $\bx_{t}^{\bu}$ is Feller continuous and strongly Markov. In particular, the measure-valued trajectory $\mathbb{P}_{\bx_{t}^{\bu}}^{\mu(t)}$ comprises of absolutely continuous measures w.r.t. Lebesgue measure.

The objective in problem (\ref{ReflectedOCPM}) is to perform the minimum control effort steering of the given initial state PDF $\rho_{0}$ at $t=0$ to the given terminal state PDF $\rho_{1}$ at $t=1$ subject to the controlled sample path dynamics (\ref{ReflectedOCPdynM}). In other words, the data of the problem consists of the domain $\overline{\mathcal{X}}$, the prior dynamics data $\bm{f},\theta$, and the two endpoint PDFs $\rho_{0},\rho_{1}$. \par
Formally, we can transcribe (\ref{ReflectedOCPM}) into the following variational problem \cite{benamou2000computational}:
\begin{subequations}  \label{SBproblem}
\begin{align}
& \underset{(\rho,\bu)\in \mathcal{P}_2(\overline{\mathcal{X}}) \times \mathcal{U}} {\text{inf}}
& & \int_{0}^{1}\!\!\!\int_{\overline{\mathcal{X}}}\frac{1}{2} \| \bu(t,\bx_{t}^{u}) \|_{2}^{2}\: \rho(t,\bx_{t}^{u})  \: \differential\bx_{t}^{u}\differential t \label{SBproblemConstraint1}\\
& \text{subject to}
& & \frac{\partial \rho }{\partial t} +\nabla \cdot \left(\rho (\bu+\bbf) \right) = \theta \Delta \rho ,  \label{SBproblemConstraint2} \\  
& & &\langle -(\bu+\bbf) \rho + \theta \nabla \rho,\bn\rangle \big |_{\partial \mathcal{X}}=0, \label{SBproblemConstraint3}\\ 
& & & \rho(0,\bx_{t}^{u}) = \rho_{0}, \quad \rho(1,\bx_{t}^{u}) = \rho_{1},  \label{SBproblemConstraint4}
\end{align}
\end{subequations}
where a PDF-valued curve $\rho(t,\cdot)\in\mathcal{P}_2(\overline{\mathcal{X}})$ if for each $t\in [0,1]$, the PDF $\rho$ is supported on $\overline{\mathcal{X}}$, and has finite second moment. In this paper, we will not focus on the rather technical direction of establishing the existence of minimizer for (\ref{SBproblem}), which can be pursued along the lines of \cite[p. 243--245]{villani2003topics}. Instead, we will formally derive the conditions of optimality, convert them to the so-called Schr\"{o}dinger system, and argue the existence-uniqueness of solutions for the same.

\subsection{Necessary Conditions of Optimality}\label{SubsecCondOptimality}
The following result summarizes how the optimal pair $(\ropt,\uopt)$ for problem (\ref{SBproblem}) can be obtained.
\begin{theorem} [\textbf{Optimal control and optimal state PDF}] 
A pair $(\ropt,\uopt)$ solving the variational problem (\ref{SBproblem}) must satisfy the system of coupled nonlinear PDEs:
\begin{subequations} \label{FPKHJB}
\begin{align}
&\frac{\partial \ropt }{\partial t} +\nabla \cdot\left(\ropt (\nabla \psi+\bbf) \right) = \theta\Delta \ropt, \label{FPK}\\ 
& \frac{\partial \psi}{\partial t} + \frac{1}{2}\lVert  \nabla \psi \rVert_{2}^2 + \langle \nabla \psi ,\bbf \rangle = -\theta\Delta \psi,\label{HJB}
\end{align}
\end{subequations}
where 
\begin{align} \label{optcntrl}
\uopt(t,\cdot) = \nabla \psi(t,\cdot),
\end{align}
subject to the boundary conditions 
\begin{subequations}
\begin{align}
    &\langle \nabla \psi, \bn \rangle \big|_{\partial \mathcal{X}}= 0, \quad \text{for all}\quad t\in[0,1], \label{HJBneumann} \\ 
    &\ropt(0,\cdot) = \rho_0, \quad \ropt(1,\cdot)= \rho_1, \label{BndryCond} \\ 
    &\langle \ropt(\nabla \psi + \bbf) -\theta \nabla \ropt,\bn \rangle \big |_{\partial \mathcal{X}}\!=0, \quad\text{for all}\quad t \in [0,1]. \label{ZeroFluxBndryCond}
\end{align}
\end{subequations}
\end{theorem}
\vspace*{0.1in}
The PDE (\ref{FPK}) is a controlled Fokker-Planck-Kolmogorov (FPK) equation, and (\ref{HJB}) is a Hamilton-Jacobi-Bellman (HJB) equation. Because the equations (\ref{FPK})-(\ref{HJB}) have one way coupling, and the boundary conditions (\ref{HJBneumann})-(\ref{ZeroFluxBndryCond}) are atypical, solving (\ref{FPKHJB}) is a challenging task in general. In the following, we show that it is possible to transform the \emph{coupled nonlinear} system (\ref{FPKHJB}) into a boundary coupled \emph{linear} system of PDEs which we refer to as the \emph{Schr\"{o}dinger system}. We will see that the resulting system paves way to a computational pipeline for solving the density steering problem with path constraints.

\subsection{Schr\"{o}dinger System}\label{SubsecSchrodingerSystem}
We now apply the Hopf-Cole transform \cite{cole1951quasi,hopf1950partial} to the system of nonlinear PDEs  (\ref{FPKHJB}).
\begin{theorem}[\textbf{Schr\"{o}dinger system}]\label{SchrodingerSystemTheorem}
Given the data $\overline{\mathcal{X}},\bbf,\theta,\rho_0,\rho_1$ for problem (\ref{SBproblem}), consider the Hopf-Cole transform $(\ropt,\psi) \mapsto (\varphi,\hat{\varphi})$ given by 
\begin{subequations}
\begin{align}
    \varphi(t,\cdot) &:= \exp\left(\psi(t,\cdot)/2\theta \right), \label{HCback}\\ 
    \hat{\varphi}(t,\cdot) &:= \ropt(t,\cdot) \exp \left(-\psi(t,\cdot)/2\theta \right), \label{HCforw}
\end{align}
\label{DefHC}
\end{subequations}
applied to (\ref{FPKHJB}) where $t\in[0,1]$. For $k\in\{0,1\}$, introduce the notation $\varphi_{k}:=\varphi(k,\cdot)$, $\hat{\varphi}_{k}:=\hat{\varphi}(k,\cdot)$. Then the pair $(\varphi,\hat{\varphi})$ satisfies the system of linear PDEs
\begin{subequations} \label{ForwardBackwardSystem}
\begin{align}
    \frac{\partial \varphi}{\partial t} &=-\langle \nabla \varphi,\bbf \rangle -\theta \Delta \varphi,  \label{BackwardEqn}\\ 
    \frac{\partial \hat{\varphi}}{\partial t} & = -\nabla \cdot( \bbf \hat{\varphi})+\theta \Delta \hat{\varphi},\label{ForwardEqn}
\end{align}
\end{subequations}
subject to the boundary conditions 
\begin{subequations} 
\begin{align}
&\varphi_{0} \hat{\varphi}_{0} = \rho_0, \quad \varphi_{1} \hat{\varphi}_{1} = \rho_1,\label{HCbndrycond1}\\
 &\langle \nabla \varphi, \bn \rangle \big|_{\partial \mathcal{X}} =
 \langle \bbf \hat{\varphi} - \theta \nabla \hat{\varphi},\bn\rangle \big|_{\partial \mathcal{X}} = 0.\label{HCbndrycond2}
 \end{align}
 \label{HCbndrycond}
\end{subequations}
For all $t\in[0,1]$, the pair $(\ropt,\uopt)$ can be recovered as 
\begin{align}
\!\ropt(t,\cdot) = \varphi(t,\cdot)\hat{\varphi}(t,\cdot), \!\!\quad \!\!\uopt(t,\cdot) = 2\theta\nabla\log\varphi(t,\cdot).
\label{recoverroptuopt}	
\end{align}
\end{theorem}

\begin{remark}
From (\ref{DefHC}), both $\varphi,\hat{\varphi}$ are nonnegative by definition, and strictly positive if $\psi$ is bounded and $\ropt$ is positive.
\end{remark}

\begin{remark} Under the regularity assumptions on $\bbf$ and $\overline{\mathcal{X}}$ stated in Section \ref{SubsecFormulation}, the process $\bx_t$ satisfying the \emph{uncontrolled} reflected It\^{o} SDE
\begin{align} \label{uncontrolledSDE}
\!\!\differential\bx_t = \bbf(t,\bx_t) \:  \differential t +\sqrt{2\theta} \: \differential \bw_t +\bm{n}(\bx_t) \: \differential \gamma_t, \: t\in[0,1],
\end{align}
is a Feller continuous strongly Markov process. Therefore, the theory of semigroups applies and the transition density of (\ref{uncontrolledSDE}) satisfies Kolmogorov's equations. Notice that the transition density or Green's function will depend on the domain $\overline{\mathcal{X}}$. In particular, we point out that (\ref{BackwardEqn}) is the \emph{backward Kolmogorov equation} in
unkonwn $\varphi$ with the corresponding \emph{Neumann boundary condition} $\langle\nabla\varphi,\bn\rangle\big |_{\partial\mathcal{X}}=0$ in (\ref{HCbndrycond2}). On the other hand, (\ref{ForwardEqn}) is the \emph{forward Kolmogorov equation} in unkonwn $\hat{\varphi}$ with the corresponding \emph{Robin boundary condition} $\langle \bbf \hat{\varphi} - \theta \nabla \hat{\varphi},\bn\rangle \big|_{\partial \mathcal{X}} = 0$ in (\ref{HCbndrycond2}). These ``backward Kolmogorov with Neumann" and ``forward Kolmogorov with Robin" system of PDE boundary value problems are coupled via the atypical boundary conditions (\ref{HCbndrycond1}). 
\end{remark}

\begin{figure}[t]
\centering
    \includegraphics[width=0.95\linewidth]{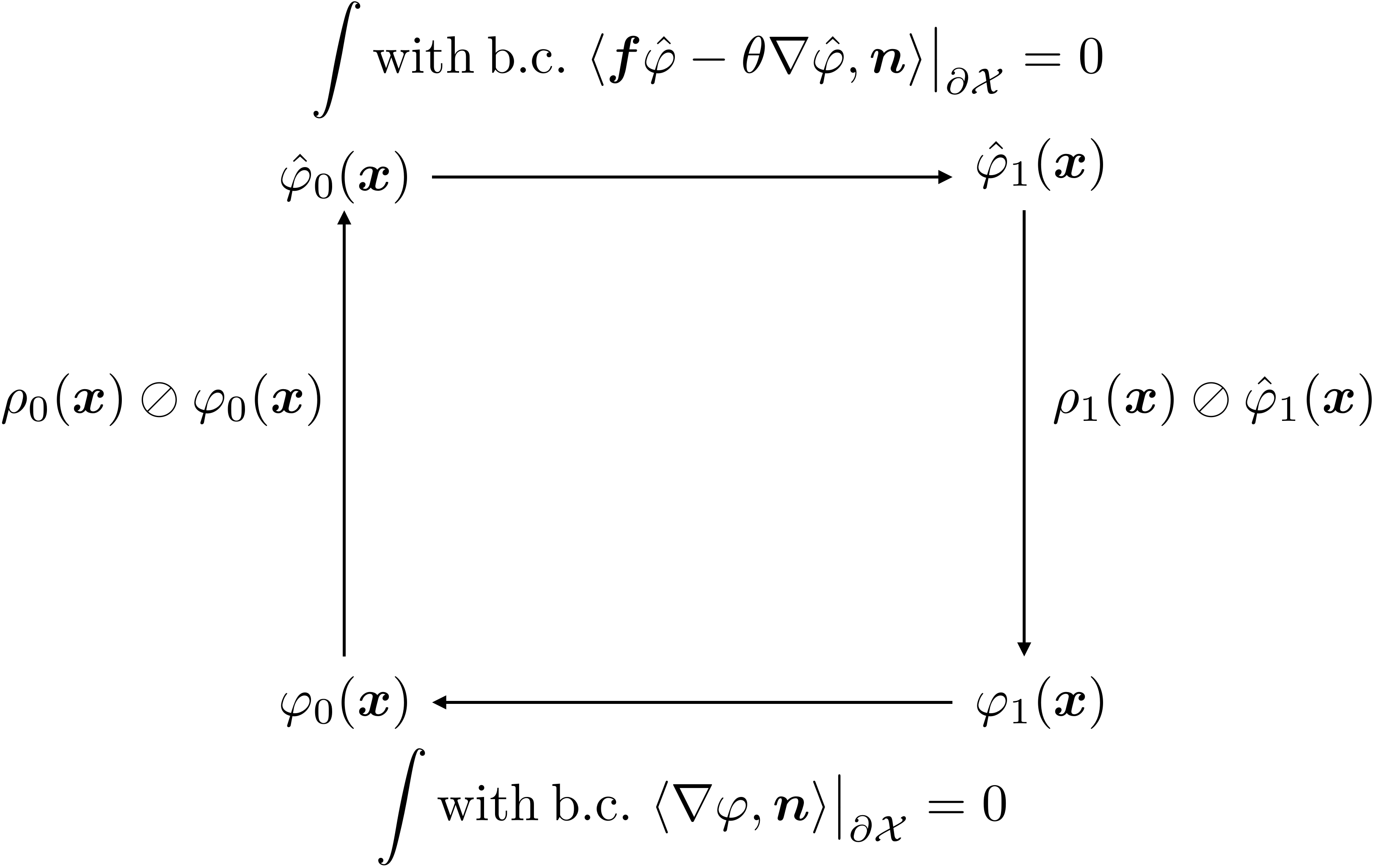}
     \caption{{\small{Schematic of the fixed point recursion for the Schr\"{o}dinger system (\ref{ForwardBackwardSystem})-(\ref{HCbndrycond}). The abbreviation ``b.c." stands for boundary condition, the symbol $\oslash$ denotes the Hadamard division.}}}
    \label{FigCommutativeDiagm}
    \vspace*{-0.23in}
\end{figure}

Theorem \ref{SchrodingerSystemTheorem} reduces finding the optimal pair $\left(\ropt,\uopt\right)$ for the RSBP to that of finding the pair\footnote{We refer to $\varphi(t,\bx_t),\hat{\varphi}(t,\bx_t)$ as the \emph{Schr\"{o}dinger factors}.} $(\varphi(t,\bx_t),\hat{\varphi}(t,\bx_t))$ associated with the uncontrolled SDE (\ref{uncontrolledSDE}). To do so, we need to compute the terminal-initial condition pair $(\varphi_1,\hat{\varphi}_0)$, which can be obtained by first making an initial guess for $(\varphi_1,\hat{\varphi}_0)$ and then performing time update by integrating the system (\ref{ForwardBackwardSystem})-(\ref{HCbndrycond2}). Using (\ref{HCbndrycond1}), this then sets up a fixed point recursion over the pair $(\varphi_1,\hat{\varphi}_0)$ (see Fig. \ref{FigCommutativeDiagm}). If this recursion converges to a unique pair, then the converged pair $(\varphi_1,\hat{\varphi}_0)$ can be used to compute the transient factors $(\varphi(t,\bx_t),\hat{\varphi}(t,\bx_t))$, and we can recover $(\ropt,\uopt)$ via (\ref{recoverroptuopt}). This computational pipeline will be pursued in this paper. 

Since the PDEs in (\ref{ForwardBackwardSystem}) are linear, and the boundary couplings in (\ref{HCbndrycond1}) are in product form, the nonnegative function pair $(\varphi_1,\hat{\varphi}_0)$ can only be unique in the projective sense, i.e., if $(\varphi_1,\hat{\varphi}_0)$ is a solution then so is $(\alpha\varphi_1,\hat{\varphi}_0/\alpha)$ for any $\alpha>0$. In \cite{chen2016entropic}, it was shown that the aforesaid fixed point recursion is in fact \emph{contractive} on a suitable cone in Hilbert's projective metric, and hence guaranteed to converge to a unique pair $(\varphi_1,\hat{\varphi}_0)$, provided that the transition density for (\ref{uncontrolledSDE}) is positive and continuous\footnote{Under the regularity assumptions on $\bbf$ and $\overline{\mathcal{X}}$ stated in Section \ref{SubsecFormulation}, the transition density for (\ref{uncontrolledSDE}) indeed satisfies
these conditions.} on $\overline{\mathcal{X}}\times\overline{\mathcal{X}}$ for all $t\in[0,1]$, and $\rho_{0},\rho_{1}$ are supported on compact subsets of $\overline{\mathcal{X}}$.

\begin{figure}[t]
\centering
    \includegraphics[width=0.9\linewidth]{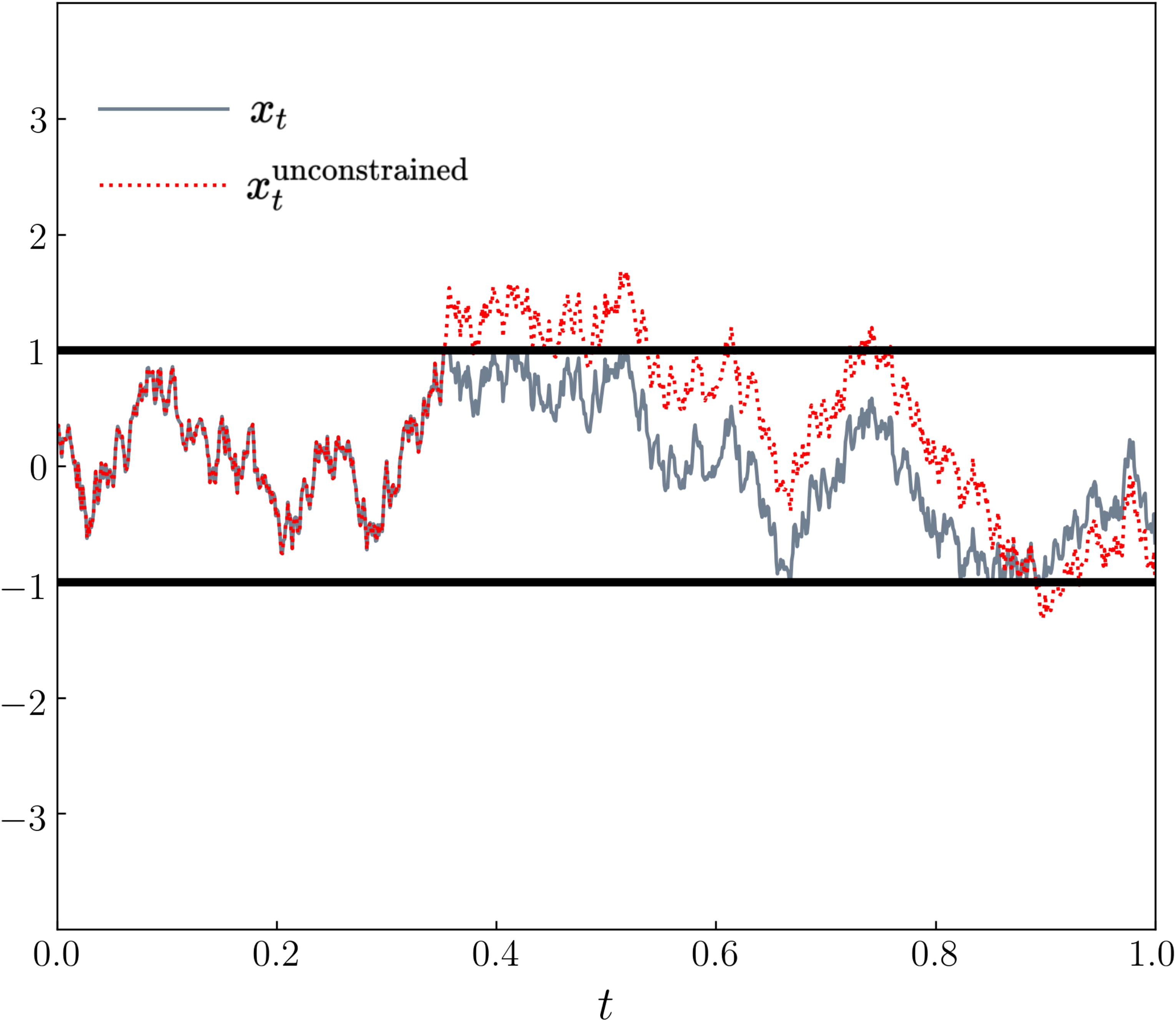}
     \caption{{\small{For $t\in[0,1]$, the \emph{solid line} shows a sample path $x_{t}$ for (\ref{TwosidedRBM}) with $[a,b]\equiv[-1,1]$, $\theta=0.5$. The \emph{dotted line} shows the corresponding unconstrained sample path $x_{t}^{\text{unconstrained}}$, computed using the two-sided Skorokhod map \cite{kruk2007explicit}.}}}
    \label{FigPathPlot}
    \vspace*{-0.23in}
\end{figure}

\section{Case Study: RSBP in 1D without Prior Drift} \label{RSBPInterval}
To illustrate the ideas presented thus far, we now consider a simple instance of problem (\ref{SBproblem}) over the state space $\overline{\mathcal{X}} =[a,b] \subset \mathbb{R}$, and with the prior drift $f\equiv 0$. That is to say, we consider the finite horizon density steering subject to the controlled two-sided reflected Brownian motion. Using some properties of the associated Markov kernel, we will show that the Schr\"{o}dinger system (\ref{ForwardBackwardSystem})-(\ref{HCbndrycond}) corresponding to this particular RSBP has a unique solution which can be obtained by the kind of fixed point recursion mentioned toward the end of Section \ref{SubsecSchrodingerSystem}. 

In this case, the Schr\"{o}dinger system (\ref{ForwardBackwardSystem})-(\ref{HCbndrycond}) reduces to
\begin{subequations} \label{ForwardBackwardHeat} 
\begin{align} 
&\frac{\partial \varphi}{\partial t} = -\theta \frac{\partial^2 \varphi}{\partial x^2}, \label{BackwardHeatPDE}\\
 &\frac{\partial \hat{\varphi}}{\partial t} = \theta \frac{\partial^2 \hat{\varphi}}{\partial x^2}, \label{ForwardHeatPDE}\\ 
 &\varphi_0\hat{\varphi}_0 = \rho_0, \quad \varphi_1\hat{\varphi}_1 = \rho_1, \label{ForwardBackwardHeatbc2} \\
&\frac{\partial \varphi}{\partial x}\bigg|_{x=a,b} = \frac{\partial \hat{\varphi}}{\partial x}\bigg|_{x=a,b} =0. \label{gradBCforwardbackwardHeat}
\end{align}
\end{subequations}
Notice that (\ref{BackwardHeatPDE})-(\ref{ForwardHeatPDE}) are the backward and forward heat PDEs, respectively, which subject to (\ref{gradBCforwardbackwardHeat}), have solutions 
\begin{subequations} \label{FormalSolution}
\begin{align}
    \varphi(x,t) &=\int_{[a ,b]} K_{\theta}(x,y,1-t) \varphi_1(y)\: \differential y  ,\quad t\leq 1, \\ 
    \hat{\varphi}(x,t) &=\int_{[a,b]} K_{\theta}(y,x,t) \hat{\varphi}_0(y)\: \differential y  , \qquad \: \: \: t \geq 0,
\end{align}
\end{subequations}
where 
\begin{align} \label{TransitionDensity}
    K_{\theta}(x,y,t)& :=\frac{1}{b-a} + \frac{2}{b-a}\sum_{m=1}^{\infty} \exp\left(-\frac{\theta \pi^2 m^2 }{(b-a)^2} t \right)  \nonumber \\
    & \times \cos\left(\frac{m \pi (x-a) }{b-a} \right)\cos\left(\frac{m \pi (y-a)}{b-a} \right)
\end{align}
is the Markov kernel or transition density \cite[Sec. 4.1]{linetsky2005transition},\cite[p. 410-411]{bhattacharya2009stochastic} associated with the uncontrolled reflected SDE 
\begin{align}
    \differential x_t = \sqrt{2\theta}\: \differential w_t + \differential L_t - \differential U_t, \quad t\in[0,1].
    \label{TwosidedRBM}
\end{align}
In (\ref{TwosidedRBM}), $L_t,U_t$ are the two local time stochastic processes \cite{glynn2018rate,harrison2013brownian} at the lower and upper boundaries respectively, which restrict $x_t$ to the interval $[a,b]$; see Fig. \ref{FigPathPlot}.
 
Combining (\ref{FormalSolution}) and (\ref{ForwardBackwardHeatbc2}), we get a system of coupled nonlinear integral equations in unknowns $(\varphi_1,\hat{\varphi}_0)$, given by 
\begin{subequations} \label{IntegralEqns}
\begin{align}
\rho_0(x) &= \hat{\varphi}_0(x) \int_{[a,b]} K_{\theta}(x,y,1) \varphi_1(y) \: \differential y, \\ 
\rho_1(x) &= \varphi_1(x) \int_{[a,b]} K_{\theta}(y,x,1) \hat{\varphi}_0(y) \: \differential y. 
\end{align}
\end{subequations} 
Clearly, solving (\ref{IntegralEqns}) is equivalent to solving (\ref{ForwardBackwardHeat}). The pair $(\varphi_1,\hat{\varphi}_0)$ can be solved from (\ref{IntegralEqns}) iteratively as a fixed point recursion with guaranteed convergence established through contraction mapping in Hilbert's projective metric;  see \cite{chen2016entropic}. The Lemma \ref{lemma} stated next will be used in the Proposition \ref{PropRBMcontraction} that follows, showing the existence and uniqueness of the pair $(\varphi_1,\hat{\varphi}_0)$ in (\ref{IntegralEqns}) as well as the fact that the aforesaid fixed point recursion is guaranteed to converge to that pair.

 \begin{lemma} 
For $0<\theta$, $a<b$, consider the transition probability density $K_{\theta}(x,y,t)$ in (\ref{TransitionDensity}). Then, 
\begin{enumerate}[i.]
    \item[(i)] $K_{\theta}(x,y,t=1)$ is continuous on the set $[a,b] \times [a,b]$. 
    \item[(ii)] $K_{\theta}(x,y,t=1)>0$ for all $(x,y) \in [a,b] \times [a,b]$.
\end{enumerate}
\label{lemma}
\end{lemma}
\begin{proposition} Given $0<\theta$, $a<b$, and the endpoint PDFs $\rho_0,\rho_1$ having compact supports $\subseteq [a,b]$. There exists a unique pair $(\varphi_1,\hat{\varphi}_0)$ that solves (\ref{IntegralEqns}), and equivalently (\ref{ForwardBackwardHeat}). Moreover, this unique pair can be computed by the fixed point recursion shown in Fig. \ref{FigCommutativeDiagm}. 
\label{PropRBMcontraction}
\end{proposition}

\begin{figure}[t]
    \centering
    \includegraphics[width=0.9\linewidth]{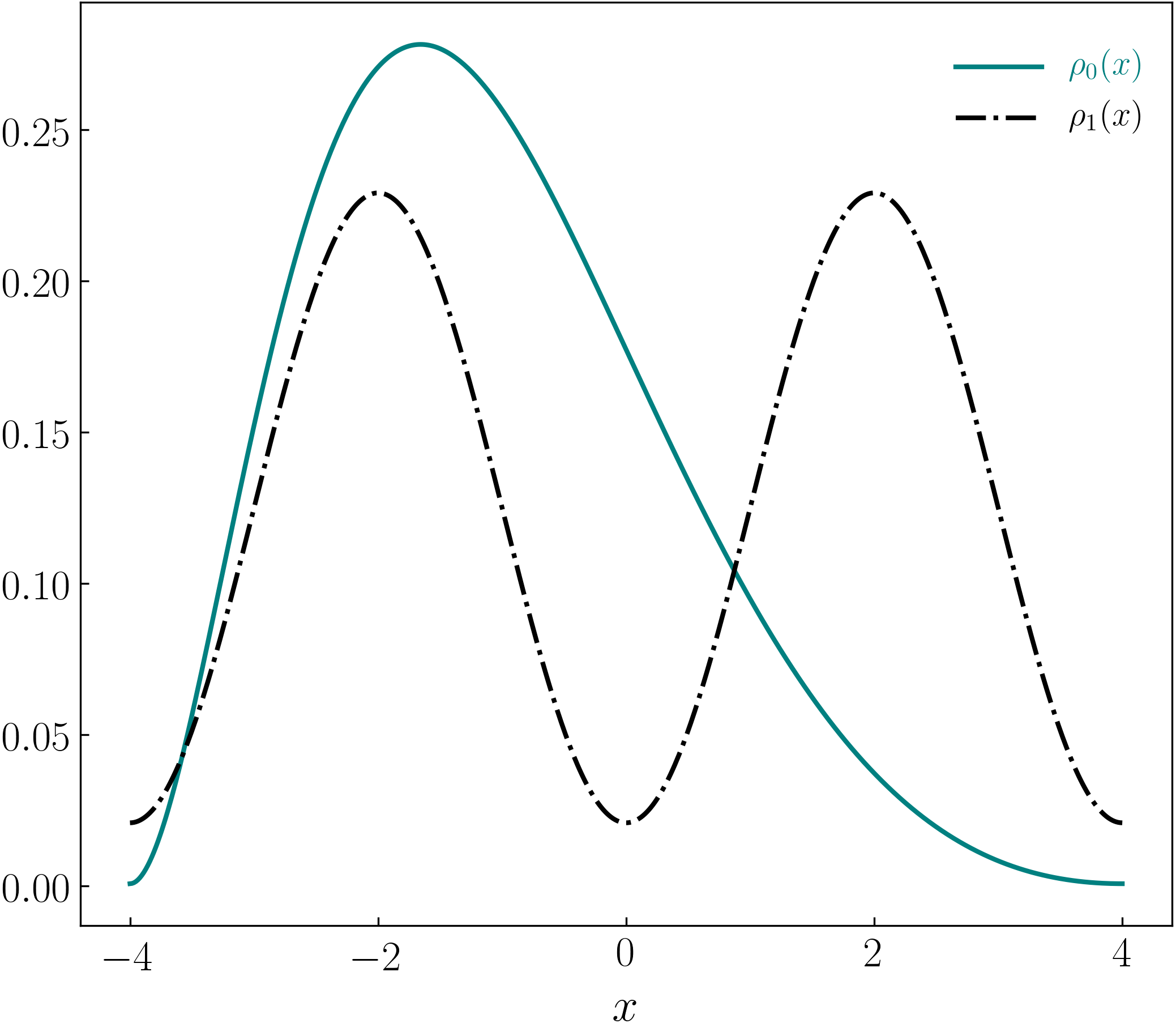}
     \caption{{\small{The endpoint PDFs $\rho_{0},\rho_{1}$ shown above are supported on $[-4,4]$, and are given by (\ref{EndpointPDFs1Dexample}).}}}
     \vspace*{-0.23in}
    \label{FigEndpointPDFs1d}
\end{figure}

\begin{figure}[htpb]
    \centering
    \includegraphics[width=0.95\linewidth]{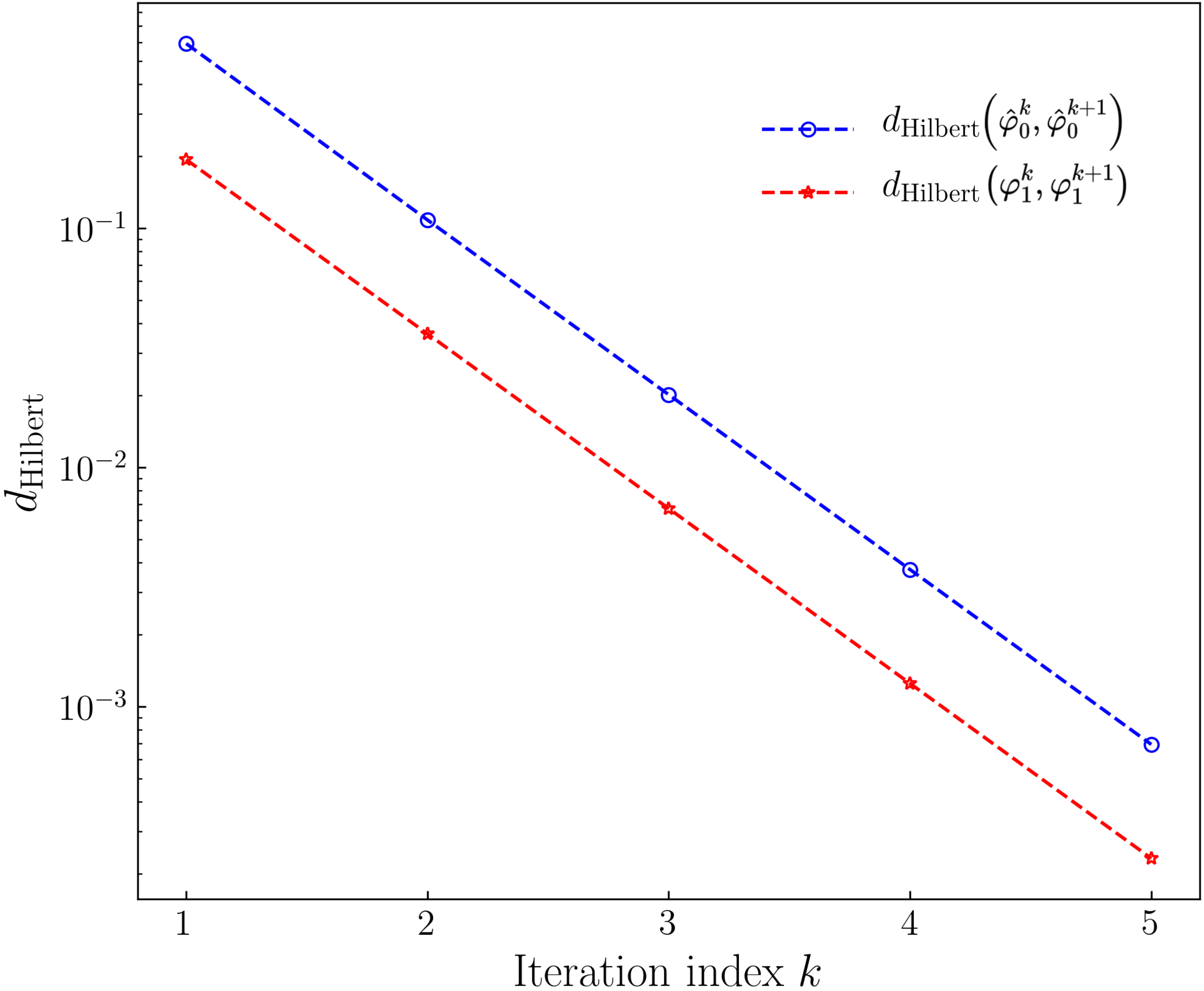}
     \caption{{\small{Convergence of the fixed point recursion over $(\varphi_{1},\hat{\varphi}_{0})$ in Hilbert's projective metric $d_{\text{Hilbert}}$.}}}
     \vspace*{-0.23in}
    \label{FigDistanceContraction}
\end{figure}

\begin{figure}[h]
    \centering 
    \includegraphics[width=0.95\linewidth]{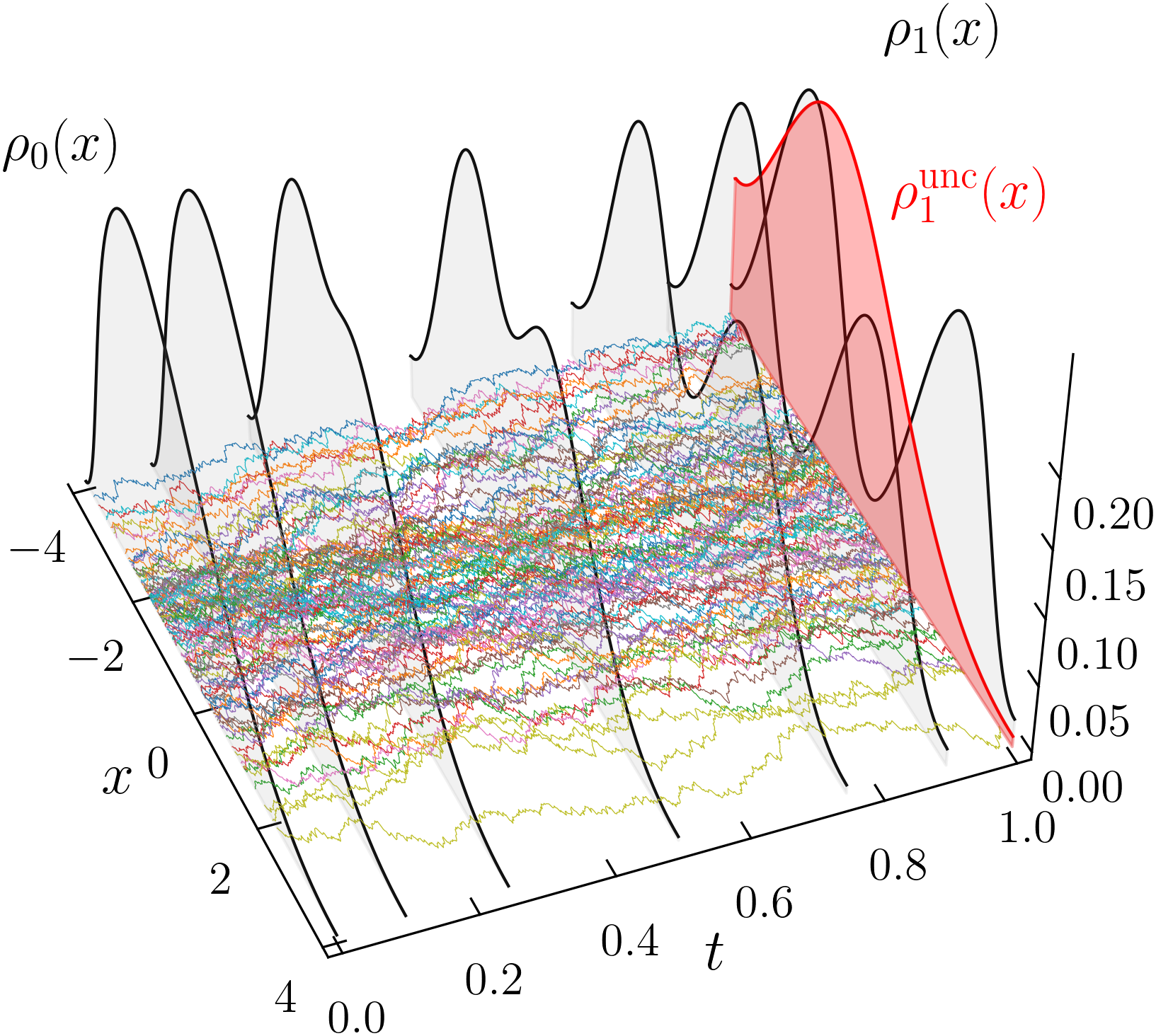}
    \caption{{\small{Shown as the black curves are the optimal controlled transient joint state PDFs $\ropt(t,x_{t}^{u})$ for steering the two-sided reflecting Brownian motion with endpoint PDFs $\rho_0,\rho_1$ as in Fig. \ref{FigEndpointPDFs1d}. The red curve $\rho_{1}^{\rm{unc}}$ is the uncontrolled state PDF at $t=1$, i.e., obtained by setting $u\equiv 0$. Also depicted are the 100 sample paths of the optimally controlled (i.e., closed-loop) reflected SDE. This simulation corresponds to the RSBP (\ref{ReflectedOCPM}) with problem data $f\equiv 0, [a,b]=[-4,4], \theta = 0.5$, and $\rho_0,\rho_1$ given by (\ref{EndpointPDFs1Dexample}).}}} 
    \vspace*{-0.2in}
    \label{FigRBM}
\end{figure}

To illustrate how the above results can be used for practical computation, consider solving the RSBP (\ref{ReflectedOCPM}) with $f\equiv 0$, $\theta=0.5$, $\overline{\mathcal{X}} = [a,b]\equiv [-4,4]$, and $\rho_{0},\rho_{1}$ as (see Fig. \ref{FigEndpointPDFs1d})
\begin{subequations}\label{EndpointPDFs1Dexample}
\begin{align}
\rho_0(x) &\propto 1 + (x^2 - 16)^{2}\exp(-x/2), \label{1dexamplerho0}\\ 
\rho_1(x) &\propto 1.2 -\cos(\pi (x+4)/2),\label{1dexamplerho1}
\end{align}
\end{subequations}
where the supports of (\ref{EndpointPDFs1Dexample}) are restricted to $[-4,4]$, and the proportionality constants are determined accordingly. For state feedback synthesis enabling this unimodal to bimodal steering over $t\in[0,1]$, we performed the fixed point recursion over the pair $\left(\varphi_{1},\hat{\varphi}_{0}\right)$ using (\ref{IntegralEqns}) with $\rho_{0},\rho_{1}$ as in (\ref{EndpointPDFs1Dexample}), and $K_{\theta}$ given by (\ref{TransitionDensity}). For numerical implementation, we truncated the infinite sum in (\ref{TransitionDensity}) after the first 100 terms. Fig. \ref{FigDistanceContraction} shows the convergence of this fixed point recursion w.r.t. Hilbert's projective metric. The converged pair $\left(\varphi_{1},\hat{\varphi}_{0}\right)$ is used to compute the transient Schr\"{o}dinger factors $(\varphi(t,\bx_t),\hat{\varphi}(t,\bx_t))$ via (\ref{FormalSolution}), and then the pair $(\ropt(t,\bx_t^u),\uopt(t,\bx_t^u))$ via (\ref{recoverroptuopt}). Fig. \ref{FigRBM} depicts the evolution of the optimal controlled transient joint state PDFs $\ropt(t,x_{t}^{u})$ as well as 100 sample paths $x_{t}^{u}$ of the optimal closed-loop reflected SDE. These sample paths were computed by applying the Euler-Maruyama scheme with time-step size $10^{-3}$. Notice from Fig. \ref{FigRBM} that (i) the closed-loop sample paths satisfy $-4\leq x_{t}^{u} \leq 4$ for all $t\in[0,1]$, and (ii) in the absence of feedback, the terminal constraint $\rho(1,x_{1}^{u})=\rho_{1}$ (given by (\ref{1dexamplerho1})) cannot be satisfied.

\section{RSBP with Prior Drift}\label{RSBPprox}
For generic $\bbf$, $\overline{\mathcal{X}}$, there is no closed-form expression of the Markov Kernel associated with (\ref{ForwardBackwardSystem})-(\ref{HCbndrycond2}). Hence, unlike the situation in Section \ref{RSBPInterval}, we cannot explicitly set up coupled integral equations of the form (\ref{IntegralEqns}), thus preventing the numerical implementation of the fixed point recursion (Fig. \ref{FigCommutativeDiagm}) via direct matrix-vector recursion. In this Section, we will show that if $\bbf$ is gradient of a potential, then we can reformulate (\ref{ForwardBackwardSystem})-(\ref{HCbndrycond}) in a way that leads to a variational recursion which in turn enables us to implement the fixed point recursion (Fig. \ref{FigCommutativeDiagm}) in an implicit manner.

\subsection{Reformulation of the Schr\"{o}dinger System}
Let $\bbf$ be a  gradient vector field, i.e., $\bbf=-\nabla V$ for some potential $V\in C^2(\overline{\mathcal{X}})$. The associated Schr\"{o}dinger system (\ref{ForwardBackwardSystem})-(\ref{HCbndrycond}) becomes
\begin{subequations}\label{ForwardBackwardGrad}
\begin{align}
    &\frac{\partial \varphi}{\partial t} =\langle \nabla \varphi,\nabla V \rangle -\theta \Delta \varphi, \label{BackwardEqnGrad} \\
    &\frac{\partial \hat{\varphi}}{\partial t} = \nabla \cdot( \nabla V \hat{\varphi})+\theta \Delta \hat{\varphi}, \label{ForwardEqnGrad} \\
    &\varphi_0\hat{\varphi}_0 = \rho_0, \quad \varphi_1\hat{\varphi}_1 = \rho_1, \\ 
 &\langle \nabla \varphi, \bn \rangle \big|_{\partial \mathcal{X}}=\langle \nabla V\hat{\varphi}+ \theta \nabla \hat{\varphi},\bn\rangle \big|_{\partial \mathcal{X}} = 0. 
 \end{align}
\end{subequations}
The idea now is to exploit the structural nonlinearities in (\ref{ForwardBackwardGrad}) to design an algorithm that allows computing the Schr\"{o}dinger factors $(\varphi,\hat{\varphi})$. To that end, the following is a crucial step.
\begin{theorem}\label{ThmReformulateSS}
Given $V\in C^2(\overline{\mathcal{X}}),\theta>0$, and $t\in[0,1]$, consider $\varphi(t,\bx_t)$ in (\ref{ForwardBackwardGrad}). Let $s:=1-t$, and define the mappings $\varphi \mapsto q \mapsto p$ given by $q(s,\bx_s):= \varphi(t,\bx_t)=\varphi(1-s,\bx_{1-s})$, $p(s,\bx_s):=q(s,\bx_s)\exp(-V(\bx_s)/\theta)$. Then $p(s,\bx_s)$ solves the PDE initial boundary value problem: 
\begin{subequations}
\begin{align}
    &\frac{\partial p}{\partial s} = \nabla \cdot (p \nabla V) + \theta \Delta p, \label{pPDE}\\  
    &p(0,\bx)= \varphi_1(\bx)\exp(-V(\bx)/\theta), \label{pInitialCondition}\\
    &\langle \nabla V p + \theta\nabla p,\bn \rangle \big|_{\partial \mathcal{X}} =0.\label{pBoundaryCondition}
\end{align}
\end{subequations}
\end{theorem}
Thanks to Theorem \ref{ThmReformulateSS}, solving (\ref{ForwardBackwardGrad}) is equivalent to solving
\begin{subequations}\label{ForwardBackwardGradREform}
\begin{align}
    &\frac{\partial p}{\partial s} = \nabla \cdot (p \nabla V) + \theta \Delta p, \label{BackwardEqnGradReform} \\
    &\frac{\partial \hat{\varphi}}{\partial t} = \nabla \cdot( \nabla V \hat{\varphi})+\theta \Delta \hat{\varphi},\label{ForwardEqnGradReform} \\
    &p(s=1,\bx)\exp(V(\bx)/\theta) \hat{\varphi}_0(\bx) = \rho_0,\nonumber\\
    &p(s=0,\bx) \exp(V(\bx)/\theta)\hat{\varphi}_1(\bx) = \rho_1, \\ 
 &\langle \nabla V p+\theta\nabla p,\bn \rangle \big|_{\partial \mathcal{X}}=\langle \nabla V\hat{\varphi}+ \theta \nabla \hat{\varphi},\bn\rangle \big|_{\partial \mathcal{X}} = 0. \label{RobinBC}
 \end{align}
\end{subequations} 
From (\ref{BackwardEqnGradReform})-(\ref{ForwardEqnGradReform}), $\varphi$ and $p$ satisfy the exact same FPK PDE with different initial conditions and integrated in different time coordinates $t$ and $s$. From (\ref{RobinBC}), $\varphi$ and $p$ satisfy the same Robin boundary condition. Therefore, a single FPK initial boundary value problem solver can be used to set up the fixed point recursion to solve for $(p_1,\hat{\varphi}_0)$, and hence $(p(s,\bx_s),\hat{\varphi}(t,\bx_t))$. From $p$, we can recover $\varphi$ as
\[\varphi(t,\bx_t)= \varphi(1-s,\bx_{1-s}) = p(s,\bx_s) \exp(-V(\bx_s)/\theta).\]

\subsection{Computation via Wasserstein Proximal Recursion} \label{proxrecur}
Building on our previous works \cite{caluya2019proximal,caluya2019gradient,caluya2019wasserstein}, we propose proximal recursions to numerically time march the solutions of the PDE initial boundary value problems (\ref{ForwardBackwardGradREform}) by exploiting certain infinite dimensional gradient descent structure. This enables us to perform the computation associated with the horizontal arrows in Fig. \ref{FigCommutativeDiagm}, and hence the fixed point recursions to solve for $(p,\hat{\varphi})$, and consequently for $(\varphi,\hat{\varphi})$. We give here a brief outline of the ideas behind these proximal recursions.

It is well-known \cite{jordan1998variational,santambrogio2017euclidean} that the flows generated by (\ref{BackwardEqnGradReform}),(\ref{ForwardEqnGradReform}),(\ref{RobinBC}) can be viewed as the gradient descent of the Lyapunov functional 
\begin{align}
F(\varrho):= \int_{\overline{\mathcal{X}}} V(\bx)\varrho(\bx) \:\differential \bx\: + \theta\int_{\overline{\mathcal{X}}} \varrho(\bx) \log \varrho(\bx) \:\differential \bx
\end{align}
w.r.t. the distance metric $W$ referred to as the (quadratic) Wassertein metric \cite{villani2003topics} on $\mathcal{P}_2(\overline{\mathcal{X}})$.
For chosen time-steps $\tau,\sigma$, this allows us to set up a variational recursion over the discrete time pair $(t_{k-1},s_{k-1}):=((k-1)\tau,(k-1)\sigma)$ as 
\begin{align}\label{ProxRecursions}
\begin{pmatrix}
    \hat{\phi}_{t_{k}} \\ 
    \varpi_{s_{k}}
\end{pmatrix} =
\begin{pmatrix}
{\rm{prox}}_{\tau F}^{W^2}(\hat{\phi}_{t_{k-1}})\\
{\rm{prox}}_{\sigma F}^{W^2}(\varpi_{s_{k-1}})
\end{pmatrix}, \quad  k \in \mathbb{N},
\end{align}
wherein the Wasserstein proximal operator %${\rm{prox}}_{h F}^{W^2}\left(\cdot\right)$ is given by
\begin{align}
{\rm{prox}}_{h F}^{W^2}(\cdot):= \underset{\varrho \in \mathcal{P}_{2}(\overline{\mathcal{X}})}{\arg\inf} \frac{1}{2}W^2(\cdot,\varrho) + h F(\varrho), \quad h>0.
\end{align}
The sequence of functions generated by the proximal recursions (\ref{ProxRecursions}) approximate the flows $(p(s,\bx_s),\hat{\varphi}(t,\bx_t))$ for (\ref{BackwardEqnGradReform}),(\ref{ForwardEqnGradReform}),(\ref{RobinBC}) in the small time step limit, i.e.,
\begin{align*}
   \hat{\phi}_{t_{k-1}} \rightarrow \hat{\varphi}(t=(k-1)\tau,\bx_{t}) \quad &\text{in $L^{1}(\overline{\mathcal{X}})$ as $\tau \downarrow 0$}, \\ 
   \varpi_{s_{k-1}} \rightarrow p(s=(k-1)\sigma,\bx_s) \quad &\text{in $L^{1}(\overline{\mathcal{X}})$ as $\sigma \downarrow 0$}.
\end{align*}
In the numerical example provided next, we solved (\ref{ProxRecursions}) using the algorithm developed in  \cite{caluya2019gradient}. 
%In other words, the variational recursion approximate the flow of the PDE as the time step goes to zero.

\subsection{Numerical Example}
\begin{figure*}[t]
\hspace{.1cm}
    \centering
    \includegraphics[width=\linewidth]{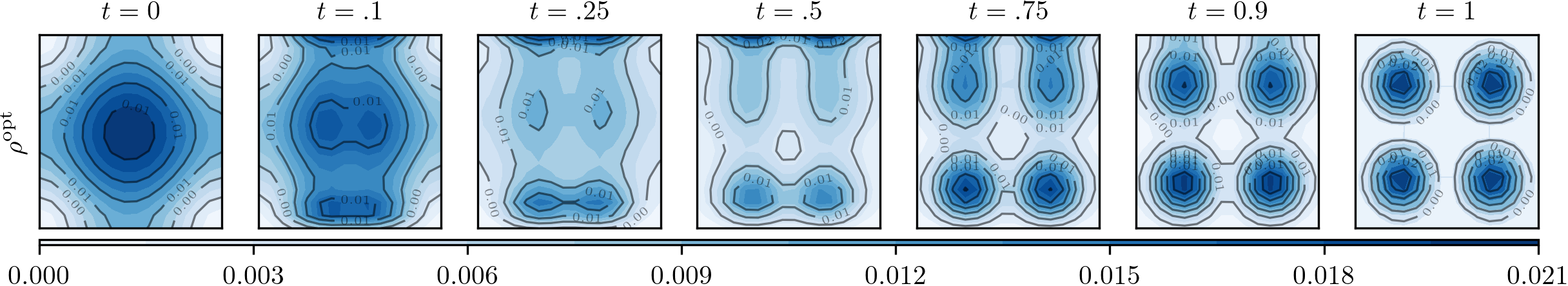}
     \caption{For the RSBP in Section \ref{proxrecur}, shown here are the contour plots of the optimal controlled joint state PDFs $\ropt(t,\bx_t^{\bu})$ over $\overline{\mathcal{X}}=[-4,4]^2$. Each subplot corresponds to a different snapshot of $\ropt$ in time. The color denotes the joint PDF value; see colorbar (dark hue = high, light hue = low).}
     \vspace*{-0.1in}
    \label{rho_con}
\end{figure*}

 We consider an instance of the RSBP with $\overline{\mathcal{X}}=[-4,4]^{2}$, $\bbf = -\nabla V$, $V(x_1,x_2):=(x_1^2+x_2^3)/5$. For 
 \begin{subequations}\label{EndpointPDFs2Dexample}
\begin{align}
\rho_0(x_1,x_2) &\propto \prod_{i=1,2} \left(1 + (x_i^2 - 16)^{2}\exp(-x_i/2)\right), \label{2dexamplerho0}\\ 
\rho_1(x_1,x_2) &\propto \prod_{i=1,2} \left(1.2 -\cos(\pi (x_i+4)/2)\right),\label{2dexamplerho1}
\end{align}
\end{subequations}
the optimal controlled joint state PDFs $\ropt(t,\bx_t^{\bu})$ are shown in Fig. \ref{rho_con}. The corresponding uncontrolled joint state PDFs $\rho^{\text{unc}}(t,\bx_t)$ are shown in Fig. \ref{rho_unc}. These results were obtained by solving (\ref{ProxRecursions}) via \cite[Sec. III.B]{caluya2019gradient} with $\tau=\sigma=10^{-3}$ to perform the fixed point recursion (Fig. \ref{FigCommutativeDiagm}) applied to (\ref{ForwardBackwardGradREform}).

\begin{figure*}[t]
\centering
\hspace{.2cm}
    \includegraphics[width=\linewidth]{rho_uncon.png}
     \caption{For the RSBP in Section \ref{proxrecur}, shown here are the contour plots of the uncontrolled joint state PDFs $\rho^{{\rm{unc}}}(t,\bx_t)$ over $\overline{\mathcal{X}}=[-4,4]^2$ starting from (\ref{2dexamplerho0}). Each subplot corresponds to a different snapshot of $\rho^{{\rm{unc}}}$ in time. The color denotes the joint PDF value; see colorbar (dark hue = high, light hue = low).}
     \vspace*{-0.15in}
    \label{rho_unc}
\end{figure*}

%%%%%%%%%%%%%%%%%%%%%%%%%%%%%%%%%%%%%%%%%%%%%%%%%%%%%%%%%%%%%%%%%%%%%%%%%%%%%%%%%%%%%%%%%%%%%%%%%%%%%%%%%%%%%%%%%%%%%

\section{Conclusions}
In this paper, we introduced the Reflected Schr\"{o}dinger Bridge Problem (RSBP) -- a stochastic optimal control problem for minimum energy feedback steering of a given joint PDF to another over finite horizon subject to reflecting boundary conditions on the controlled state trajectories. Combining our prior work on Wasserstein proximal recursions with some recent results on contraction mapping associated with the Schr\"{o}dinger system, we provide a computational pipeline for optimal feedback synthesis. Numerical examples are given to highlight the proposed framework.

%%%%%%%%%%%%%%%%%%%%%%%%%%%%%%%%%%%%%%%%%%%%%%%%%%%%%%%%%%%%%%%%%%%%%%%%%%%%%%%%%%%%%%%%%%%%%%%%%%%%%%%%%%%%%%%%%%%%%

\appendix
\subsection{Proof of Theorem 1}
\noindent The necessary conditions for optimality (\ref{FPKHJB}) can be deduced using the Lagrange multiplier theorem in Banach spaces; see \cite[Ch.4.14, Proposition 1]{zeidler1995applied}. This theorem allows us set up an augmented Lagrangian associated with (\ref{SBproblem}) and perform pointwise minimization to derive (\ref{FPKHJB}). 

To apply this in our context, define the function spaces
\begin{align*}
\mathcal{P}_{01} &:=\{\rho(t,\cdot)\in\mathcal{P}_{2}\left(\overline{\mathcal{X}}\right) \mid \rho(0,\cdot)=\rho_0, \rho(1,\cdot)=\rho_1\},\\
\widetilde{\mathcal{P}}_{01} &:= \mathcal{P}_{01} \cap L^{2}(H^{1}([0,1];\overline{\mathcal{X}})) \cap \dot{H}^{1}\left(\left(H^{1}([0,1];\overline{\mathcal{X}})\right)^{*}\right),\\
X &:= \widetilde{\mathcal{P}}_{01} \times  L^{2}([0,1]\times\overline{\mathcal{X}}), \quad Y:=  L^{2}(H^{-1}([0,1];\overline{\mathcal{X}})),
\end{align*}
where $[0,1]$ denotes the time interval, and $L^{2}\left(\cdot\right)$ denotes the space of square integrable functions. The notation $L^{2}([0,1];H^{1}(\overline{\mathcal{X}}))$ stands for the Sobolev space of functions having first order weak derivatives w.r.t. $\bx_t^{\bu}\in\overline{\mathcal{X}}$, and finite $L^{2}$ norms w.r.t. $t\in[0,1]$. Furthermore, $\dot{H}^{1}\left([0,1];\left(H^{1}(\overline{\mathcal{X}})\right)^{*}\right):= \{\phi(t,\cdot)\in L^{2}\left([0,1]\right) \mid \frac{\partial\phi}{\partial t} \in L^{2}\left([0,1]\right),\phi\in\left(H^{1}(\overline{\mathcal{X}})\right)^{*}\}$, wherein $\left(H^{1}(\overline{\mathcal{X}})\right)^{*}$ denotes the dual space of the Sobolev space $H^{1}(\overline{\mathcal{X}})$. We denote the dual space of $\widetilde{\mathcal{P}}_{01}$ as $\widetilde{\mathcal{P}}_{01}^{*}$. In the definition of $Y$, the notation $H^{-1}\left(\overline{\mathcal{X}}\right)$ stands for the space of all linear functionals on $H_{0}^{1}\left(\mathcal{X}\right):=\{\phi\in H^{1}\left(\mathcal{X}\right), \; \text{and vanishes on}\; \partial\mathcal{X}\}$. Then, in (\ref{SBproblemConstraint1}), the objective functional $F: X \mapsto \mathbb{R}$, and is given by 
\begin{align}
F(\rho,\bu):= \int_{\overline{\mathcal{X}}} \int_{0}^{1} \frac{1}{2} \lVert \bu(t,\bx_t^{\bu}) \rVert_2^2 \:\rho(t,\bx_t^{\bu}) \: \differential t \: \differential \bx_t^{\bu}. 
\end{align}
The constraint is a mapping  $G: X \mapsto Y$ given by
{\small{\begin{align}
&G(\rho,\bu)(\psi)\!\!:= \!\!\!\int_{\overline{\mathcal{X}}}\!\!\psi(1,\bx_t^{\bu}) \rho(1,\bx_t^{\bu})\differential \bx_t^{\bu} \!- \!\!\int_{\overline{\mathcal{X}}}\!\!\psi(0,\bx_t^{\bu}) \rho(0,\bx_t^{\bu})\differential\bx_t^{\bu} \nonumber\\&
-\!\!\int_{\overline{\mathcal{X}}} \!\int_{0}^{1}\!\frac{\partial \psi}{\partial t} \rho \: \differential \bx_t^{\bu} \: \differential t  \! +\!\! 
\int_{\overline{\mathcal{X}}}\!\int_{0}^{1}\!\!\!\psi( \nabla \cdot(\rho\bu+\bbf) - \theta \Delta \rho) \: \differential \bx_t^{\bu} \: \differential t,
\end{align}}}
where we used (\ref{SBproblemConstraint3}) so that the boundary
terms vanish in the integration by parts. Following \cite[p. 112-114]{albi2017mean}, one can show that $G'_{\rho}(\rho,\bu)$ and $G'_{\bu}(\rho,\bu)$ (where $'$ denotes derivative w.r.t. the subscripted variable) are surjective, and hence by \cite[Ch.4.14, Proposition 1]{zeidler1995applied}, there exists $\psi\in Y^{*}=L^2([0,1];H_0^1(\overline{\mathcal{X}}))$. This result allows us to perform pointwise minimization of the augmented Lagrangian
\begin{align*}
    &\mathscr{L}(\rho,\bu,\psi):=\!\!\!\int_{0}^{1}\!\!\!\int_{\mathcal{X}}\frac{1}{2} \| \bu(t,\bx_t^{\bu}) \|_{2}^{2}\:\rho(t,\bx_t^{\bu})  \: \differential\bx_t^{\bu} \differential t\:+\nonumber \\
    &\underbrace{\int_{\mathcal{X}}\!\!\int_{0}^{1}\!\!\psi \frac{\partial \rho}{\partial t} \differential t \differential\bx_t^{\bu}}_{\rm{term  1}}+ \underbrace{\int_{0}^{1}\!\!\!\int_{\mathcal{X}}\!\!\psi \left(  \nabla \cdot(\rho \bu+\bbf) - \theta \Delta \rho \right)  \: \differential\bx_t^{\bu} \differential t}_{\rm{term  2}}.
  \end{align*}
By performing integration by parts in $t$, term $1$ becomes 
\begin{align*}
\!\!\displaystyle\int_{\overline{\mathcal{X}}}\!\!\left( \underbrace{\psi(1,\bx_t^{\bu})\rho(1,\bx_t^{\bu})-\psi(0,\bx_t^{\bu})\rho(0,\bx_t^{\bu})}_{{\text{constant w.r.t. $(\rho,\bm{u})$}}} - \!\!\int_{0}^1\!\!\!\frac{\partial \psi}{\partial t} \rho\differential t \right)\!\!\differential\bx_t^{\bu}.
\end{align*}
For term 2, we perform integration by parts w.r.t $\bx_t^{\bu}$, impose the boundary condition (\ref{SBproblemConstraint3}), and thereby deduce that $\mathscr{L}$ (up to an additive constant) equals 
{\small{\begin{align} 
\label{SimplifiedLagrangian1}
\!\!\int_{0}^{1}\!\!\!\int_{\overline{\mathcal{X}}}\!\!\left( \frac{1}{2} \| \bu \|_{2}^{2} -\frac{\partial \psi}{\partial t} - \langle \nabla \psi, \bu+\bbf \rangle - \theta\Delta \psi\!\right)\!\rho \: \differential\bx_t^{\bu} \:\differential t.
\end{align}}}
Pointwise minimization of (\ref{SimplifiedLagrangian1}) w.r.t $\bu$ while fixing $\rho$, gives the optimal control (\ref{optcntrl}). Substituting (\ref{optcntrl}) back into (\ref{SimplifiedLagrangian1}) and equating the resulting expression to zero results in the dynamic programming equation 
\begin{align*} 
\!\!\int_{0}^{1}\!\!\!\int_{\overline{\mathcal{X}}}\!\!\left( -\frac{\partial \psi}{\partial t} - \frac{1}{2}\lVert  \nabla \psi\rVert^2 - \langle \nabla \psi, \bbf \rangle- \theta\Delta \psi\right)\!\!\rho\:\differential\bx_t^{\bu}\:\differential t = 0.
\end{align*}
Since the above holds for arbitrary $\rho$, we must have 
\begin{align*}
   \frac{\partial \psi}{\partial t} + \frac{1}{2}\lVert  \nabla \psi \rVert^2 + \langle \nabla \psi , \bbf \rangle +\theta\Delta \psi= 0,
\end{align*}
which is indeed the HJB PDE (\ref{HJB}). Substituting (\ref{optcntrl}) in (\ref{SBproblemConstraint1}) yields the FPK PDE (\ref{FPK}).

The Neumann boundary condition (\ref{HJBneumann}) follows directly (see \cite{lions1984optimal}). The endpoint conditions (\ref{BndryCond}) follow from (\ref{SBproblemConstraint4}). The Robin boundary condition (\ref{ZeroFluxBndryCond}) is obtained by combining (\ref{optcntrl}) with (\ref{SBproblemConstraint3}).\qed

%%%%%%%%%%%%%%%%%%%%%%%%%%%%%%%%%%%%%%%%%%%%%%%%%%%%%%%%%%%%%%%%%%%%%%%%%%%%%%%%%%%%%%%%%%%%%%%%%%%%%%%%%%%%%%%%%%%%%

\subsection{Proof of Theorem 2}
\noindent The system of linear PDEs (\ref{ForwardBackwardSystem}) are obtained via straightforward but tedious computation detailed in \cite[Appendix B]{caluya2019wasserstein}. The boundary conditions (\ref{HCbndrycond1}) follow by setting $t=0,1$ in (\ref{DefHC}). To derive (\ref{HCbndrycond2}), evaluate (\ref{HCback}) at a boundary point $\bx_{\text{bdy}} \in \partial \mathcal{X}$. In the resulting expression, take the natural log to both sides and then take the gradient w.r.t. $\bx_{\text{bdy}}$, to get
\begin{align}
    \nabla \psi(t,\bx_{\text{bdy}}) = 2\theta \frac{\nabla \varphi(t,\bx_{\text{bdy}})}{\varphi(t,\bx_{\text{bdy}})}. 
\label{gradpsibdy}    
\end{align}
In both sides of (\ref{gradpsibdy}), we take the inner product with the normal vector $\bn(\bx_{\text{bdy}})$, and use (\ref{HJBneumann}), to obtain $\langle \nabla \varphi, \bn \rangle \big|_{\partial \mathcal{X}} = 0$, as in (\ref{HCbndrycond2}). To deduce the second equality in (\ref{HCbndrycond2}), we evaluate (\ref{HCforw}) at $\bx_{\text{bdy}}$, and then as before, take the natural log followed by the gradient w.r.t. $\bx_{\text{bdy}}$, and invoke (\ref{HJBneumann}) to arrive at
\begin{align}
\frac{\langle\nabla\hat{\varphi}(t,\bx_{\text{bdy}}),\bn(\bx_{\text{bdy}})\rangle}{\hat{\varphi}(t,\bx_{\text{bdy}})} = \frac{\langle\ropt(t,\bx_{\text{bdy}}),\bn(\bx_{\text{bdy}})\rangle}{\ropt(t,\bx_{\text{bdy}})}.
\label{Intermed9bsecondequality}	
\end{align}
Using (\ref{HJBneumann}) again in (\ref{ZeroFluxBndryCond}), the right-hand-side of (\ref{Intermed9bsecondequality}) simplifies to $\langle\bbf(t,\bx_{\text{bdy}}),\bn(\bx_{\text{bdy}})\rangle/\theta$, and thus yields the second equality in (\ref{HCbndrycond2}). Finally, (\ref{recoverroptuopt}) follows from (\ref{DefHC}) and (\ref{optcntrl}).\qed

%%%%%%%%%%%%%%%%%%%%%%%%%%%%%%%%%%%%%%%%%%%%%%%%%%%%%%%%%%%%%%%%%%%%%%%%%%%%%%%%%%%%%%%%%%%%%%%%%%%%%%%%%%%%%%%%%%%%%

\subsection{Proof of Lemma 1}
\noindent \emph{Proof of (i):} To demonstrate continuity, it suffices to show that the infinite sum in (\ref{TransitionDensity}) for $K_{\theta}(x,y,t=1)$ converges uniformly on $[a,b]\times [a,b]$. For $k\in\mathbb{N}$, let 
\begin{align}
    f_k(x,y) := &\exp\left(-\frac{\theta \pi^2 k^2 }{(b-a)^2}  \right)  \cos\left(\frac{k \pi (x-a) }{b-a} \right) \nonumber \\ 
    & \times \cos\left(\frac{k \pi (y-a)}{b-a} \right),
\end{align}
and notice that
\begin{align} \label{WeierstrassInequality}
    \lvert f_{k}(x,y) \rvert \leq  M_k \quad \text{for all}\quad(x,y) \in [a,b] \times [a,b],
\end{align}
where $M_k := \exp\left(-\theta \pi^2 k^2/(b-a)^2\right)$. Furthermore,
\begin{align}
    \lim_{k\rightarrow \infty} \bigg \lvert \frac{M_{k+1}}{M_k}\bigg \rvert =  \lim_{k\rightarrow \infty} \exp\left(-\frac{\theta \pi^2 (2k+1)}{(b-a)^2} \right) = 0.
\end{align}
By the ratio test \cite[Ch. 3, Theorem 3.34]{rudin1976principles}, we then have
\begin{align} \label{Convergent bound}
    \sum_{k=1}^{\infty} M_k < \infty. 
\end{align}
From (\ref{WeierstrassInequality}) and (\ref{Convergent bound}), the Weierstrass M-test \cite[Ch. 7, Theorem 7.10]{rudin1976principles} implies that $\sum_{k=1}^{\infty} f_k(x,y)$ is uniformly convergent for all $(x,y)\in [a,b] \times [a,b]$, and the resulting sum must converge to a continuous function. Therefore, $K_{\theta}(x,y,t=1)$ is continuous for all $(x,y) \in [a,b] \times [a,b]$. 

\noindent\emph{Proof of (ii):}  To establish positivity, set $r:=b-a,\tx := x-a, \ty := y-a$. Using basic trigonometry and the Euler's identity, we find that
\begin{align}
&\frac{1}{r} + \frac{2}{r}\sum_{m=1}^{\infty} \exp\left(-\frac{\theta \pi^2 m^2 }{r^2}  \right)\left[\cos\left(\frac{m \pi \tx }{r} \right)\cos\left(\frac{m \pi \ty}{r} \right)\right] \nonumber\\
    % &= \frac{1}{r} + \frac{1}{r} \sum_{m=1}^{\infty}\exp\left(-\frac{\theta \pi^2 m^2 }{r^2}  \right) \nonumber \\
    % & \qquad\qquad\times \left[\cos\left( \frac{ m\pi(\tx+\ty) }{r}\right) +\cos\left( \frac{ m\pi(\tx-\ty) }{r}\right) \right] \nonumber \\ 
    % &=  \frac{1}{r} + \frac{1}{2r} \sum_{m=1}^{\infty}\exp\left(-\frac{\theta \pi^2 m^2 }{r^2}  \right) \nonumber\\
    % &\qquad\times \left[\exp{\left(\frac{im \pi (\tx+\ty)}{r} \right)}+\exp{\left(\frac{-im \pi (\tx+\ty)}{r} \right)} \right. \left. \right. \nonumber \\ 
    % &\qquad+ \exp{\left(\frac{im \pi (\tx-\ty)}{r} \right)} +\left.\exp{\left(\frac{-im \pi (\tx-\ty)}{r} \right)}\right] \nonumber \\ 
    % &=\frac{1}{r} + \frac{1}{2r}  \sum_{\substack{m=-\infty,\\ m\neq 0}}^{\infty}\exp\left(-\frac{\theta \pi^2 m^2 }{r^2}  \right)   \nonumber \\
    % &\qquad\times  \left[\exp\left({\frac{im\pi(\tx+\ty)}{r}}\right) \exp\left({\frac{im\pi (\tx-\ty)}{r}}\right) \right]\nonumber \\ 
     &=  \sum_{m=-\infty}^{\infty} \frac{1}{2r}\exp\left(-\frac{\theta \pi^2 m^2 }{r^2}  \right)\exp\left({\frac{im\pi(\tx+\ty)}{r}}\right)+ \nonumber \\  
     &\: \: \: \:  \sum_{m=-\infty}^{\infty} \frac{1}{2r} \exp\left(-\frac{\theta \pi^2 m^2 }{r^2}  \right) \exp\left({\frac{im\pi(\tx-\ty)}{r}}\right).  
\label{sumfinalexpression}
\end{align}
Let 
{\small{\begin{align}
    g(m) := \underbrace{\frac{1}{2r}\exp\left(-\frac{\theta \pi^2 m^2 }{r^2}  \right)}_{=: g_1(m)} \underbrace{\exp \left(  \frac{im\pi(\tx+\ty)}{r} \right)}_{=: g_2(m)},
\end{align}}}
and denote the Fourier transforms of $g_1(m),g_2(m)$ as $\widehat{g}_1(\widehat{m}),\widehat{g}_2(\widehat{m})$, respectively. Notice that 
\begin{align}
    \widehat{g}_1(\widehat{m}) &= \frac{1}{\sqrt{4 \pi \theta}} \exp \left(  -\frac{r^2 \widehat{m}^2}{\theta }\right),\\
     \widehat{g}_2(\widehat{m}) &= \delta\left(\widehat{m} - \frac{\tx+\ty}{2r}\right),
\end{align}
where $\delta(\cdot)$ denotes the Dirac delta, and hence by the convolution theorem, the Fourier transform of $g$ is
{\small{\begin{align}
    \widehat{g}(\widehat{m}) &= \frac{1}{\sqrt{4\pi\theta}} \exp\left(-\frac{r^2}{\theta} \left(\widehat{m}-\frac{(\tx+\ty)}{2r}\right)^2\right) \nonumber \\ 
    &= \frac{1}{\sqrt{4\pi\theta}} \exp\left(-\frac{(2\widehat{m}r-\tx-\ty)^2}{4\theta } \right).
\end{align}}}
Invoking the Poisson summation formula \cite[Ch. 4, Theorem 2.4]{stein2010complex}, we deduce
{\small{\begin{align}
    \sum_{m=-\infty}^{\infty} g(m) =  \sum_{\widehat{m}=-\infty}^{\infty} \widehat{g}(\widehat{m}),
\end{align}}}
implying that the infinite sum in (\ref{sumfinalexpression}) is equal to
{\small{\begin{align}
\frac{1}{\sqrt{4\pi\theta}} \sum_{\widehat{m}=-\infty}^{\infty} &\left[\exp\left(-\frac{(2\widehat{m}r-\tx-\ty)^2}{4\theta}\right) + \right.  \nonumber \\ 
 & \:   \left. \exp \left(-\frac{(2\widehat{m}r-\tx+\ty)^2}{4\theta } \right) \right],
\end{align}}}
which is obviously positive. Therefore, $K_{\theta}(x,y,t=1)$ is positive for all $(x,y) \in [a,b] \times [a,b]$.\qed

%%%%%%%%%%%%%%%%%%%%%%%%%%%%%%%%%%%%%%%%%%%%%%%%%%%%%%%%%%%%%%%%%%%%%%%%%%%%%%%%%%%%%%%%%%%%%%%%%%%%%%%%%%%%%%%%%%%%%

\subsection{Proof of Proposition 1}
\noindent Using Lemma \ref{lemma} and that $[a,b]$ is a compact metric space, the hypotheses of \cite[Proposition 4 and Theorem 8]{chen2016entropic}  are satisfied. Therefore, the solution pair $(\varphi_1,\varphi_0)$ exists and is unique in the projective sense. Furthermore, the fixed point recursion being contractive in Hilbert's projective metric, converges (by contraction mapping theorem) to this pair.\qed

%%%%%%%%%%%%%%%%%%%%%%%%%%%%%%%%%%%%%%%%%%%%%%%%%%%%%%%%%%%%%%%%%%%%%%%%%%%%%%%%%%%%%%%%%%%%%%%%%%%%%%%%%%%%%%%%%%%%%

\subsection{Proof of Theorem 3}
\noindent The derivation of (\ref{pPDE})-(\ref{pInitialCondition}) follows \cite[Appendix C]{caluya2019wasserstein}. By substituting the identity $\nabla q = \exp(V/\theta)(\nabla p + p\nabla V/\theta)$ in $0 = \langle\nabla\varphi,\bn\rangle = \langle\nabla q,\bn\rangle$, we find (\ref{pBoundaryCondition}).\qed

\bibliographystyle{IEEEtran}
\bibliography{references.bib}

% Generated by IEEEtran.bst, version: 1.14 (2015/08/26)
\begin{thebibliography}{10}
\providecommand{\url}[1]{#1}
\csname url@samestyle\endcsname
\providecommand{\newblock}{\relax}
\providecommand{\bibinfo}[2]{#2}
\providecommand{\BIBentrySTDinterwordspacing}{\spaceskip=0pt\relax}
\providecommand{\BIBentryALTinterwordstretchfactor}{4}
\providecommand{\BIBentryALTinterwordspacing}{\spaceskip=\fontdimen2\font plus
\BIBentryALTinterwordstretchfactor\fontdimen3\font minus
  \fontdimen4\font\relax}
\providecommand{\BIBforeignlanguage}[2]{{%
\expandafter\ifx\csname l@#1\endcsname\relax
\typeout{** WARNING: IEEEtran.bst: No hyphenation pattern has been}%
\typeout{** loaded for the language `#1'. Using the pattern for}%
\typeout{** the default language instead.}%
\else
\language=\csname l@#1\endcsname
\fi
#2}}
\providecommand{\BIBdecl}{\relax}
\BIBdecl

\bibitem{hotz1987covariance}
A.~Hotz and R.~E. Skelton, ``Covariance control theory,'' \emph{International
  Journal of Control}, vol.~46, no.~1, pp. 13--32, 1987.

\bibitem{skelton1989covariance}
R.~Skelton and M.~Ikeda, ``Covariance controllers for linear continuous-time
  systems,'' \emph{International Journal of Control}, vol.~49, no.~5, pp.
  1773--1785, 1989.

\bibitem{grigoriadis1997minimum}
K.~M. Grigoriadis and R.~E. Skelton, ``Minimum-energy covariance controllers,''
  \emph{Automatica}, vol.~33, no.~4, pp. 569--578, 1997.

\bibitem{chen2015optimal1}
Y.~Chen, T.~T. Georgiou, and M.~Pavon, ``Optimal steering of a linear
  stochastic system to a final probability distribution, {Part I},'' \emph{IEEE
  Transactions on Automatic Control}, vol.~61, no.~5, pp. 1158--1169, 2015.

\bibitem{chen2015optimal2}
------, ``Optimal steering of a linear stochastic system to a final probability
  distribution, {Part II},'' \emph{IEEE Transactions on Automatic Control},
  vol.~61, no.~5, pp. 1170--1180, 2015.

\bibitem{chen2016optimal}
------, ``Optimal transport over a linear dynamical system,'' \emph{IEEE
  Transactions on Automatic Control}, vol.~62, no.~5, pp. 2137--2152, 2016.

\bibitem{villani2003topics}
C.~Villani, \emph{Topics in optimal transportation}.\hskip 1em plus 0.5em minus
  0.4em\relax American Mathematical Soc., 2003, no.~58.

\bibitem{schrodinger1931umkehrung}
E.~Schr{\"o}dinger, ``{\"U}ber die umkehrung der naturgesetze,''
  \emph{Sitzungsberichte der Preuss. Phys. Math. Klasse}, vol.~10, pp.
  144--153, 1931.

\bibitem{schrodinger1932theorie}
------, ``Sur la th{\'e}orie relativiste de l'{\'e}lectron et
  l'interpr{\'e}tation de la m{\'e}canique quantique,'' in \emph{Annales de
  l'institut Henri Poincar{\'e}}, vol.~2, no.~4, 1932, pp. 269--310.

\bibitem{halder2016finite}
A.~Halder and E.~D. Wendel, ``Finite horizon linear quadratic {G}aussian
  density regulator with {W}asserstein terminal cost,'' in \emph{2016 American
  Control Conference (ACC)}.\hskip 1em plus 0.5em minus 0.4em\relax IEEE, 2016,
  pp. 7249--7254.

\bibitem{bakolas2018finite}
E.~Bakolas, ``Finite-horizon covariance control for discrete-time stochastic
  linear systems subject to input constraints,'' \emph{Automatica}, vol.~91,
  pp. 61--68, 2018.

\bibitem{okamoto2019input}
K.~Okamoto and P.~Tsiotras, ``Input hard constrained optimal covariance
  steering,'' \emph{arXiv preprint arXiv:1903.10964, IEEE Conference on
  Decision and Control, Nice, France}, 2019.

\bibitem{bakolas2017covariance}
E.~Bakolas, ``Covariance control for discrete-time stochastic linear systems
  with incomplete state information,'' in \emph{2017 American Control
  Conference (ACC)}.\hskip 1em plus 0.5em minus 0.4em\relax IEEE, 2017, pp.
  432--437.

\bibitem{caluya2019finite}
K.~F. Caluya and A.~Halder, ``Finite horizon density control for static state
  feedback linearizable systems,'' \emph{arXiv preprint arXiv:1904.02272},
  2019.

\bibitem{caluya2019finiteMulti}
------, ``Finite horizon density steering for multi-input state feedback
  linearizable systems,'' \emph{arXiv preprint arXiv:1909.12511, to appear in
  2020 American Control Conference}, 2019.

\bibitem{caluya2019wasserstein}
------, ``Wasserstein proximal algorithms for the {S}chr{\"o}dinger bridge
  problem: Density control with nonlinear drift,'' \emph{arXiv preprint
  arXiv:1912.01244}, 2019.

\bibitem{okamoto2018optimal}
K.~Okamoto, M.~Goldshtein, and P.~Tsiotras, ``Optimal covariance control for
  stochastic systems under chance constraints,'' \emph{IEEE Control Systems
  Letters}, vol.~2, no.~2, pp. 266--271, 2018.

\bibitem{ikeda1961construction}
N.~Ikeda, ``On the construction of two-dimensional diffusion processes
  satisfying {W}entzell’s boundary conditions and its application to boundary
  value problems,'' \emph{Memoirs of the College of Science, University of
  Kyoto. Series A: Mathematics}, vol.~33, no.~3, pp. 367--427, 1961.

\bibitem{watanabe1971stochastic}
S.~Watanabe, ``On stochastic differential equations for multi-dimensional
  diffusion processes with boundary conditions,'' \emph{Journal of Mathematics
  of Kyoto University}, vol.~11, no.~1, pp. 169--180, 1971.

\bibitem{lions1984stochastic}
P.-L. Lions and A.-S. Sznitman, ``Stochastic differential equations with
  reflecting boundary conditions,'' \emph{Communications on Pure and Applied
  Mathematics}, vol.~37, no.~4, pp. 511--537, 1984.

\bibitem{harrison1987multidimensional}
J.~M. Harrison and R.~J. Williams, ``Multidimensional reflected {B}rownian
  motions having exponential stationary distributions,'' \emph{The Annals of
  Probability}, pp. 115--137, 1987.

\bibitem{glynn2018rate}
P.~W. Glynn and R.~J. Wang, ``On the rate of convergence to equilibrium for
  reflected {B}rownian motion,'' \emph{Queueing Systems}, vol.~89, no. 1-2, pp.
  165--197, 2018.

\bibitem{harrison2013brownian}
J.~M. Harrison, \emph{Brownian models of performance and control}.\hskip 1em
  plus 0.5em minus 0.4em\relax Cambridge University Press, 2013.

\bibitem{pilipenko2014introduction}
A.~Pilipenko, \emph{An introduction to stochastic differential equations with
  reflection}.\hskip 1em plus 0.5em minus 0.4em\relax Universit{\"a}tsverlag
  Potsdam, 2014, vol.~1.

\bibitem{skorokhod1961stochastic}
A.~V. Skorokhod, ``Stochastic equations for diffusion processes in a bounded
  region,'' \emph{Theory of Probability \& Its Applications}, vol.~6, no.~3,
  pp. 264--274, 1961.

\bibitem{skorokhod1962stochastic}
------, ``Stochastic equations for diffusion processes in a bounded region.
  {II},'' \emph{Theory of Probability \& Its Applications}, vol.~7, no.~1, pp.
  3--23, 1962.

\bibitem{kruk2007explicit}
L.~Kruk, J.~Lehoczky, K.~Ramanan, S.~Shreve \emph{et~al.}, ``An explicit
  formula for the {S}korokhod map on [0, a],'' \emph{The Annals of
  Probability}, vol.~35, no.~5, pp. 1740--1768, 2007.

\bibitem{stroock1971diffusion}
D.~W. Stroock and S.~S. Varadhan, ``Diffusion processes with boundary
  conditions,'' \emph{Communications on Pure and Applied Mathematics}, vol.~24,
  no.~2, pp. 147--225, 1971.

\bibitem{benamou2000computational}
J.-D. Benamou and Y.~Brenier, ``A computational fluid mechanics solution to the
  {M}onge-{K}antorovich mass transfer problem,'' \emph{Numerische Mathematik},
  vol.~84, no.~3, pp. 375--393, 2000.

\bibitem{cole1951quasi}
J.~D. Cole, ``On a quasi-linear parabolic equation occurring in aerodynamics,''
  \emph{Quarterly of Applied Mathematics}, vol.~9, no.~3, pp. 225--236, 1951.

\bibitem{hopf1950partial}
E.~Hopf, ``The partial differential equation $u_{t} + uu_{x} = \mu_{xx}$,''
  \emph{Communications on Pure and Applied mathematics}, vol.~3, no.~3, pp.
  201--230, 1950.

\bibitem{chen2016entropic}
Y.~Chen, T.~Georgiou, and M.~Pavon, ``Entropic and displacement interpolation:
  a computational approach using the {H}ilbert metric,'' \emph{SIAM Journal on
  Applied Mathematics}, vol.~76, no.~6, pp. 2375--2396, 2016.

\bibitem{linetsky2005transition}
V.~Linetsky, ``On the transition densities for reflected diffusions,''
  \emph{Advances in Applied Probability}, vol.~37, no.~2, pp. 435--460, 2005.

\bibitem{bhattacharya2009stochastic}
R.~N. Bhattacharya and E.~C. Waymire, \emph{Stochastic processes with
  applications}.\hskip 1em plus 0.5em minus 0.4em\relax SIAM, 2009, vol.~61.

\bibitem{caluya2019proximal}
K.~F. Caluya and A.~Halder, ``Proximal recursion for solving the
  {F}okker-{P}lanck equation,'' in \emph{2019 American Control Conference
  (ACC)}.\hskip 1em plus 0.5em minus 0.4em\relax IEEE, 2019, pp. 4098--4103.

\bibitem{caluya2019gradient}
K.~Caluya and A.~Halder, ``Gradient flow algorithms for density propagation in
  stochastic systems,'' \emph{IEEE Transactions on Automatic Control}, 2019.

\bibitem{jordan1998variational}
R.~Jordan, D.~Kinderlehrer, and F.~Otto, ``The variational formulation of the
  {F}okker-{P}lanck equation,'' \emph{SIAM journal on mathematical analysis},
  vol.~29, no.~1, pp. 1--17, 1998.

\bibitem{santambrogio2017euclidean}
F.~Santambrogio, ``$\{$Euclidean, metric, and Wasserstein$\}$ gradient flows:
  an overview,'' \emph{Bulletin of Mathematical Sciences}, vol.~7, no.~1, pp.
  87--154, 2017.

\bibitem{zeidler1995applied}
E.~Zeidler, \emph{Applied functional analysis: main principles and their
  applications}.\hskip 1em plus 0.5em minus 0.4em\relax Springer Science \&
  Business Media, 1995, vol. 109.

\bibitem{albi2017mean}
G.~Albi, Y.-P. Choi, M.~Fornasier, and D.~Kalise, ``Mean field control
  hierarchy,'' \emph{Applied Mathematics \& Optimization}, vol.~76, no.~1, pp.
  93--135, 2017.

\bibitem{lions1984optimal}
P.~Lions, ``Optimal control of reflected diffusion processes,'' in
  \emph{Filtering and Control of Random Processes}.\hskip 1em plus 0.5em minus
  0.4em\relax Springer, 1984, pp. 157--163.

\bibitem{rudin1976principles}
W.~Rudin, \emph{Principles of Mathematical Analysis}.\hskip 1em plus 0.5em
  minus 0.4em\relax Third Edition, McGraw-Hill, USA, 1976.

\bibitem{stein2010complex}
E.~M. Stein and R.~Shakarchi, \emph{Complex analysis}.\hskip 1em plus 0.5em
  minus 0.4em\relax Princeton University Press, 2010, vol.~2.

\end{thebibliography}

\end{document}